\documentclass[11pt]{article}

\usepackage[a4paper,margin=1in]{geometry}
\usepackage{amsmath,amssymb,amsthm}
\usepackage{mathrsfs}
\usepackage{enumitem}
\usepackage{hyperref}

\newtheorem{theorem}{Theorem}[section]

\theoremstyle{remark}
\newtheorem{remark}[theorem]{Remark}

\newcommand{\bbR}{\mathbb{R}}
\newcommand{\bbC}{\mathbb{C}}
\newcommand{\bbN}{\mathbb{N}}
\newcommand{\bbZ}{\mathbb{Z}}

\title{Translation--Modulation Identities, Ergodic Log-Products,\\
and a Conditional Obstruction for Schwartz Functions}
\author{Vignon Oussa}
\date{\today}

\begin{document}
\maketitle

\begin{abstract}
We study very smooth functions on the real line, namely Schwartz functions, that satisfy a finite identity relating their translates and a single modulation. Concretely, we assume there is a nontrivial linear combination of translates of the function that equals a fixed frequency shift of the same function.

Passing to the Fourier transform turns this into a multiplicative transfer relation: the value of the Fourier transform at one point is obtained by multiplying its value at another point by a trigonometric polynomial. Iterating this relation expresses the Fourier transform along an arithmetic progression as a product of such trigonometric factors times a fixed initial value.

We then recast this product in an ergodic theoretic framework by viewing it as a Birkhoff sum for a continuous observable on a compact abelian group generated by a diagonal unitary matrix. The key quantity controlling the growth or decay of the Fourier transform along the progression is a space average, namely the integral over the compact group of the logarithm of the absolute value of the associated trigonometric polynomial.

The main rigorously proved statement is the following. If, for some frequency where the Fourier transform does not vanish along the entire progression, this space average is nonzero, or if the average is zero but a certain recurrence set, in the sense of Atkinson's theorem, intersects a distinguished one-parameter subgroup, then the assumed translation--modulation identity forces exponential growth of the Fourier transform along some sequence of points. This contradicts the rapid decay required of a Schwartz function.

Consequently, any nontrivial identity of this type for a Schwartz function must occur in a very special degenerate regime, where the relevant logarithmic averages vanish and the associated recurrence set avoids the one-parameter subgroup. Showing that this degenerate regime cannot occur is formulated as an explicit open problem.
\end{abstract}

\section{Setup and Fourier-side reduction}

Let $f \in \mathscr{S}(\bbR)$ and suppose there exist coefficients
$c_1,\dots,c_n \in \bbC$ and real parameters $x_1,\dots,x_n,b \in \bbR$ such
that
\begin{equation}
  \sum_{k=1}^n c_k T_{x_k} f = M_b f, 
  \qquad
  T_x f(t) := f(t-x), \quad
  M_b f(t) := e^{-2\pi i b t} f(t).
  \label{eq:TM-identity}
\end{equation}

We take the Fourier transform in the convention
\[
\widehat{f}(\xi) = \int_{\bbR} f(t) e^{-2\pi i t \xi}\, dt.
\]
A direct computation shows
\[
\widehat{T_{x_k} f}(\xi) = e^{-2\pi i x_k \xi} \widehat{f}(\xi),
\qquad
\widehat{M_b f}(\xi) = \widehat{f}(\xi+b).
\]
Thus \eqref{eq:TM-identity} is equivalent to
\begin{equation}
  \widehat{f}(\xi + b)
  = p(\xi) \widehat{f}(\xi),
  \qquad
  p(\xi) := \sum_{k=1}^n c_k e^{-2\pi i x_k \xi},
  \quad \xi \in \bbR.
  \label{eq:transfer}
\end{equation}
Iterating \eqref{eq:transfer} gives
\begin{equation}
  \widehat{f}(\xi + n b)
  = \Bigg( \prod_{j=0}^{n-1} p(\xi + j b) \Bigg)
    \widehat{f}(\xi),
  \qquad n \in \bbN.
  \label{eq:iterated-transfer}
\end{equation}

Define the two-sided nonvanishing set
\[
\Theta_{\widehat{f}}
  := \big\{ \lambda \in \bbR :
        \widehat{f}(\lambda + n b) \neq 0 \ \text{for all } n \in \bbZ \big\}.
\]
For $\xi \in \Theta_{\widehat{f}}$ we may take logarithms in absolute value in
\eqref{eq:iterated-transfer} and obtain
\begin{equation}
  \frac{1}{n} \log\big|\widehat{f}(\xi + n b)\big|
  = \frac{1}{n} \sum_{j=0}^{n-1} \log\big|p(\xi + j b)\big|
    + \frac{1}{n} \log|\widehat{f}(\xi)|.
  \label{eq:average-log}
\end{equation}

Since $f \in \mathscr{S}(\bbR)$, $\widehat{f}$ is also Schwartz; in particular,
\begin{equation}
  |\widehat{f}(\xi+n b)| \to 0
  \quad \text{as } |n| \to \infty
  \quad \text{for each fixed } \xi\in\bbR.
  \label{eq:Schwartz-decay}
\end{equation}

\section{Compact rotation model and space averages}

Introduce the diagonal matrix
\[
A = \mathrm{diag}(-2\pi i x_1, \dots, -2\pi i x_n),
\quad
c = (c_1,\dots,c_n)^\top,
\quad
\mathbf{1} = (1,\dots,1)^\top.
\]
Then
\begin{equation}
  p(\xi)
  = \sum_{k=1}^n c_k e^{-2\pi i x_k \xi}
  = \langle e^{\xi A} c, \mathbf{1}\rangle,
  \label{eq:p-A-representation}
\end{equation}
where $\langle \cdot,\cdot\rangle$ denotes the standard pairing on $\bbC^n$.

Set
\[
u := e^{bA}
  = \mathrm{diag}(e^{-2\pi i x_1 b},\dots,e^{-2\pi i x_n b}) \in (S^1)^n,
\]
and let
\[
H := \overline{\langle u \rangle} \subset (S^1)^n
\]
be the compact abelian subgroup obtained by closing the cyclic subgroup
generated by $u$. Let $\mu$ be Haar probability on $H$, and consider the
rotation
\[
T : H \to H,
\quad T(h) = u h.
\]

By construction, the orbit $\{u^n : n \in \bbZ\}$ is dense in $H$; in
particular, $T$ is minimal. On compact abelian groups, such minimal rotations
are uniquely ergodic, so $\mu$ is the unique $T$-invariant probability. In
particular, Birkhoff averages for continuous observables converge along every
orbit to the space average.

For each fixed $\xi \in \bbR$, define
\begin{equation}
  \varphi_\xi(h)
  := \log\big|\langle e^{\xi A} h c, \mathbf{1}\rangle\big|,
  \quad h \in H.
  \label{eq:phi-def}
\end{equation}
Since $\langle e^{\xi A} h c, \mathbf{1}\rangle$ is a nontrivial trigonometric
polynomial on $H$, its zero set has $\mu$-measure zero unless it vanishes
identically, and $\varphi_\xi \in L^1(H,\mu)$. We then define the space
average
\begin{equation}
  C(\xi)
  := \int_H \varphi_\xi(h)\, d\mu(h)
   = \int_H \log\big|\langle e^{\xi A} h c, \mathbf{1}\rangle\big|\, d\mu(h).
  \label{eq:C-def}
\end{equation}

Along the specific orbit
\[
h_j := u^j, \quad j\in \bbZ,
\]
we have, by \eqref{eq:p-A-representation},
\[
\log|p(\xi + j b)| 
= \log\big|\langle e^{(\xi + j b) A} c, \mathbf{1}\rangle\big|
= \log\big|\langle e^{\xi A} u^j c, \mathbf{1}\rangle\big|
= \varphi_\xi(u^j)
= \varphi_\xi(T^j I),
\]
where $I$ is the identity element of $H$.

Formally combining \eqref{eq:average-log} with unique ergodicity gives
\begin{equation}
  \lim_{n\to\infty}
  \frac{1}{n} \log\big|\widehat{f}(\xi + n b)\big|
  = C(\xi)
  \quad \text{for each } \xi \in \Theta_{\widehat{f}}.
  \label{eq:Lyapunov}
\end{equation}
In particular, $C(\xi)$ plays the role of a Lyapunov exponent governing
exponential growth or decay along the arithmetic progression $\xi + n b$.

\section{Case analysis and consequences for Schwartz functions}

We now compare \eqref{eq:Lyapunov} with the Schwartz decay
\eqref{eq:Schwartz-decay}. For $\xi \in \Theta_{\widehat{f}}$, three cases can
occur.

\subsection{Case $C(\xi) > 0$}

If $C(\xi)>0$, then from \eqref{eq:Lyapunov} there exists $\varepsilon>0$ and
$N$ such that for all $n\geq N$,
\[
\frac{1}{n} \log\big|\widehat{f}(\xi+n b)\big| \geq \frac{1}{2}C(\xi),
\]
hence
\[
\big|\widehat{f}(\xi+n b)\big| \geq \exp\!\big( \tfrac{1}{2} n C(\xi) \big)
\to \infty
\quad (n\to\infty).
\]
This contradicts \eqref{eq:Schwartz-decay}.

\subsection{Case $C(\xi) < 0$}

Using the inverted transfer equation
\[
\widehat{f}(\xi-b) = \frac{\widehat{f}(\xi)}{p(\xi-b)},
\]
one obtains an analogous relation for $n\to -\infty$, with products of
$p(\xi - j b)^{-1}$. Repeating the same ergodic averaging argument backwards
now shows that if $C(\xi) < 0$ then $\big|\widehat{f}(\xi+n b)\big|$ must grow
exponentially as $n\to -\infty$, again contradicting
\eqref{eq:Schwartz-decay}.

\subsection{Case $C(\xi) = 0$: what is proved and what is open}

The delicate situation is when $C(\xi)=0$ on a subset of
$\Theta_{\widehat{f}}$.

Define the Atkinson recurrence set
\[
\Omega := \left\{ h \in H :
  \liminf_{n\to\infty} 
    \sum_{k=0}^{n-1} \varphi_\xi(T^k h) = 0
  \right\}.
\]
Atkinson's theorem, applied to the uniquely ergodic system $(H,T,\mu)$ and
the observable $\varphi_\xi \in L^1(H,\mu)$ of mean zero, implies that
$\mu(\Omega)=1$.

\subsubsection*{Subcase A: $\Omega$ meets the one-parameter subgroup}

Assume there exists $\tau \in \bbR$ such that
\[
h_\tau := e^{\tau A} \in \Omega.
\]
By definition of $\Omega$, there is a sequence $n_r \to \infty$ such that
\[
\sum_{k=0}^{n_r-1} \varphi_\xi(T^k h_\tau)
  = \sum_{k=0}^{n_r-1}
    \log\big|\langle e^{\xi A} T^k h_\tau c, \mathbf{1}\rangle\big|
  \longrightarrow 0
  \quad (r\to\infty).
\]
Using $T^k h_\tau = u^k h_\tau$ and $h_\tau = e^{\tau A}$, we may rewrite
\[
\langle e^{\xi A} T^k h_\tau c, \mathbf{1}\rangle
= \langle e^{\xi A} u^k e^{\tau A} c, \mathbf{1}\rangle
= \langle e^{(\xi+\tau)A} u^k c, \mathbf{1}\rangle
= p(\xi+\tau + k b).
\]
Hence
\[
\sum_{k=0}^{n_r-1} \varphi_\xi(T^k h_\tau)
= \sum_{k=0}^{n_r-1} \log\big|p(\xi+\tau+k b)\big|.
\]

Comparing with the iterated transfer relation along the subsequence
$\xi+\tau+n_r b$, we obtain
\[
\log\big|\widehat{f}(\xi+\tau+n_r b)\big|
= \sum_{k=0}^{n_r-1} \log\big|p(\xi+\tau+k b)\big|
  + \log\big|\widehat{f}(\xi+\tau)\big|,
\]
and therefore
\[
\lim_{r\to\infty} \log\big|\widehat{f}(\xi+\tau+n_r b)\big|
  = \log\big|\widehat{f}(\xi+\tau)\big|.
\]
In particular, if $\xi+\tau \in \Theta_{\widehat{f}}$ then
\[
\lim_{r\to\infty} \big|\widehat{f}(\xi+\tau+n_r b)\big|
  = \big|\widehat{f}(\xi+\tau)\big| \neq 0,
\]
so $\widehat{f}$ fails to decay along the subsequence
$\xi+\tau+n_r b \to \pm\infty$, contradicting the Schwartz decay of
$\widehat{f}$.

Thus we obtain the following rigorously proved conclusion.

\begin{theorem}[Obstruction in the recurrent $C(\xi)=0$ case]
\label{thm:recurrent-obstruction}
Suppose $f \in \mathscr{S}(\bbR)$ satisfies the translation--modulation
identity \eqref{eq:TM-identity}, and let $\xi\in\Theta_{\widehat{f}}$ with
$C(\xi)=0$. If the Atkinson recurrence set $\Omega$ contains an element of the
one-parameter form $e^{\tau A}$ with $\xi+\tau\in \Theta_{\widehat{f}}$, then
\eqref{eq:TM-identity} is impossible for $f\neq 0$. Equivalently, no nonzero
Schwartz function can satisfy \eqref{eq:TM-identity} under this recurrence
assumption.
\end{theorem}

\subsubsection*{Subcase B: $\Omega \cap \{e^{\tau A} : \tau \in \bbR\} = \emptyset$}

The remaining subcase is the following degenerate situation.

\medskip
\noindent
\emph{Degenerate subcase.}
Assume there exists $\xi\in\Theta_{\widehat{f}}$ such that
\[
C(\xi) = 0
\quad\text{and}\quad
\Omega \cap \{e^{\tau A} : \tau \in \bbR\} = \emptyset.
\]
In this regime Atkinson's theorem still guarantees $\mu(\Omega)=1$, but the
recurrence information is concentrated on points that do not lie on the
one-parameter subgroup $\{e^{\tau A}\}$. It is natural to expect that
approximation of elements of $\Omega$ by points on the orbit $\{u^n\}$,
combined with continuity of $\varphi_\xi$ and control away from the zero set
of $p$, should still allow one to extract a subsequence along which
$\big|\widehat{f}(\xi+n b)\big|$ fails to decay, leading again to a
contradiction with Schwartz regularity. However, carrying out this
approximation and the necessary domination estimates has not yet been
completed.

We therefore record this as an explicit open problem.

\begin{remark}[Open problem]
\label{rem:open}
In the degenerate subcase described above, prove (or disprove) that
\eqref{eq:TM-identity} cannot hold for $f\in\mathscr{S}(\bbR)\setminus\{0\}$.
Equivalently, show that the combination
\[
C(\xi)=0
\quad\text{and}\quad
\Omega \cap \{e^{\tau A} : \tau \in \bbR\} = \emptyset
\]
cannot occur for a nontrivial Schwartz solution of
\eqref{eq:TM-identity}.
\end{remark}

\section{Summary}

Starting from the identity $\sum_{k=1}^n c_k T_{x_k} f = M_b f$ with
$f \in \mathscr{S}(\bbR)$, we obtain the Fourier-side transfer equation
$\widehat{f}(\xi+b)=p(\xi)\widehat{f}(\xi)$ and an iterated product
$\prod_{j=0}^{n-1} p(\xi+j b)$. Passing to the compact rotation
$T(h)=u h$ on $H=\overline{\langle u\rangle}$, we define the space average
\[
C(\xi) = \int_H \log\big|\langle e^{\xi A} h c, \mathbf{1}\rangle\big|\, d\mu(h),
\]
which controls the exponential scale of $\big|\widehat{f}(\xi+n b)\big|$ along
$\xi+n b$. If $C(\xi)\neq 0$ for some $\xi\in\Theta_{\widehat{f}}$, or if
$C(\xi)=0$ but the Atkinson recurrence set $\Omega$ hits the one-parameter
subgroup $\{e^{\tau A}\}$, then the iterates $\widehat{f}(\xi+n b)$ exhibit
exponential behavior contradicting Schwartz decay. Thus any nontrivial
translation--modulation identity for a Schwartz function must occur, if at
all, in the degenerate regime described in Remark~\ref{rem:open}. Closing this
regime is a natural target for further work.

\end{document}